\documentclass[a4paper,12pt]{article}

\usepackage{amsmath,amsfonts,graphicx}

\numberwithin{equation}{section}

\title{The congruent number descent of Komotu, Watanabe and Wada.}

\author{Allan J. MacLeod,\\Statistics, O.R. and Mathematics Group (retired),\\
University of the West of Scotland,\\High St., Paisley,\\Scotland.  PA1 2BE\\
(e-mail: peediejenn@hotmail.com)}

\date{}

\begin{document}

\maketitle

\begin{abstract}
We present an elementary exposition of the descent method used by Komoto, Watanabe and Wada to show
$N=42553$ is a congruent number. We generalize to $N=7p$ with $p$ prime, and try to present the method
such that an amateur could understand and perform their own calculations.
\end{abstract}

\newpage

\section{Introduction}
Whilst surfing the Internet, I came across \cite{kww}. The number $N=42553=7*6079$ comes from the paper by Nemenzo\cite{nem}, where it
is stated to be the smallest unproven congruent/non-congruent number (at that time). Komato, Watanabe and Wada prove it is congruent by
explicitly computing a point on the related elliptic curve. I was fascinated by the simplicity of the method, but some details were not given.
I now describe the method, for $N=7*\mbox{prime}$, in full detail.

The elliptic curve is
\begin{equation*}
Y^2=X^3-N^2X
\end{equation*}
and with $X=du^2/v^2$, $Y=duw/v^3$ where $d,u,v,w,\in \mathbb{Z}$, we have
\begin{equation*}
d^2w^2=d^3u^4-N^2dv^4
\end{equation*}
with $d$ squarefree and $d|N$, $\gcd(u,v)=1$, $\gcd(w,v)=1$.

Suppose $N=7p$ where $p$ is prime. Pick $d=-p$ so
\begin{equation}
w^2=-p(u^4-49v^4)
\end{equation}
and, since $p|w$ we can set $w=pz$, giving
\begin{equation}
pz^2=(49v^4-u^4)=(7v^2+u^2)(7v^2-u^2)
\end{equation}

Assume $u,v$ both odd so $7v^2+u^2$ is a multiple of $8$. Thus set
\begin{equation*}
7v^2+u^2=8b^2 \hspace{2cm} 7v^2-u^2=2pa^2
\end{equation*}
giving $z=ab$, and
\begin{equation*}
7v^2=4b^2+pa^2 \hspace{2cm} u^2=4b^2-pa^2
\end{equation*}

The latter identity can be written
\begin{equation*}
pa^2=4b^2-u^2=(2b+u)(2b-u)
\end{equation*}
so set
\begin{equation*}
2b+u=pc^2 \hspace{2cm} 2b-u=d^2
\end{equation*}
and $a=c\,d$.

Thus,
\begin{equation*}
4b=pc^2+d^2 \hspace{2cm} 2u=pc^2-d^2
\end{equation*}

From $7v^2=4b^2+pa^2$, we get the identity,
\begin{equation}
p^2c^4+6pc^2d^2+d^4=28v^2
\end{equation}
which is what the method attempts to solve. Let $e=c^2, f=d^2$ so we have the quadratic
\begin{equation*}
p^2e^2+6pef+f^2=28v^2
\end{equation*}
which we can write
\begin{equation}
(e,f,v) \left(\begin{array}{lrr}p^2&3p&0\\3p&1&0\\0&0&-28\end{array} \right) \left( \begin{array}{r}e\\f\\v \end{array} \right) = 0
\end{equation}

Let $(1,f_0,v_0)$ be a solution. $e=1$ makes things a lot easier later on, though we cannot always find such a solution. If we can, $f_0$
must be odd.

Define
\begin{equation}
\left( \begin{array}{l}e\\f\\v \end{array} \right) = \left( \begin{array}{lrr}1&0&0\\f_0&1&0\\v_0&0&1 \end{array} \right) \left( \begin{array}{l}g\\h\\i \end{array} \right)
\end{equation}
so
\begin{equation*}
2f_0gh+2g(3hp-28iv_0)+h^2-28i^2=0
\end{equation*}

Thus $h$ must be even and we set $h=2j$. This gives the equation
\begin{equation}
(g,j,i) \left(\begin{array}{lrr}0&(f_0+3p)/2&-7v_0\\(f_0+3p)/2&1&0\\-7v_0&0&-7\end{array} \right) \left( \begin{array}{r}g\\j\\i \end{array} \right) = 0
\end{equation}

What we do next is to transform the matrix so that the $(1,3)$ element becomes zero and a non-zero element is in the $(2,3)$ place.

Let $A=(f_0+3p)/2$ and $B=-7v_0$, and assume $\gcd(A,B)=1$. Then there exists integers $\alpha, \beta$ such that $\alpha \, A + \beta \, B=1$.
Then, set
\begin{equation}
\left( \begin{array}{l}g\\j\\i \end{array} \right) = \left( \begin{array}{lrr}1&0&0\\0&\alpha&-B\\0&\beta&A \end{array} \right) \left( \begin{array}{l}g\\k\\l \end{array} \right)
\end{equation}

This gives
\begin{equation}
(g,k,l) \left(\begin{array}{lrr}0&1&0\\1&\alpha^2-7\beta^2&-\alpha B-7\beta A\\0&-\alpha B-7\beta A&-14p^2\end{array} \right) \left( \begin{array}{r}g\\k\\l \end{array} \right) = 0
\end{equation}
and we can simplify by setting $k=2m, l=2n$ to have
\begin{equation*}
m( g+(\alpha^2-7\beta^2)m+7(2\alpha v_0-\beta(f_0+3p))n )=14p^2n^2
\end{equation*}

Define
\begin{equation}\label{reqn}
r=g+(\alpha^2-7\beta^2)m+7(2\alpha v_0-\beta(f_0+3p))n
\end{equation}
so that our equation to solve reduces to
\begin{equation}
m\,r=14p^2\, n^2
\end{equation}

Set $m=14p^2s^2$, $r=t^2$ and $n=st$. Then
\begin{equation*}
c^2=e=g=t^2+7(\beta (f_0+3p)-2\alpha v_0)ts+14p^2(7\beta^2-\alpha^2)s^2
\end{equation*}
so
\begin{equation*}
c^2=\left( t+\frac{7(\beta (f_0+3p)-2\alpha v_0)s}{2} \right)^2+Ks^2
\end{equation*}
where $K$ is a constant to be determined.

With the help of a symbolic algebra package we find $K=-7$.

Defining
\begin{equation*}
t+\frac{7(\beta (f_0+3p)-2\alpha v_0)s}{2}=W
\end{equation*}
gives
\begin{equation*}
7s^2=W^2-c^2=(W+c)(W-c)
\end{equation*}

Then, set $W+c=56y^2$, $W-c=2x^2$ so that $s=4xy$, $W=28y^2+x^2$ and $c=28y^2-x^2$.

We also have $d^2=f=f_0g+h=f_0 g+4(\alpha m-Bn)$.

Using \eqref{reqn} and making various substitutions we finally arrive at the quartic
\begin{equation}
d^2=f_0 x^4+112v_0x^3y +
\end{equation}
\begin{equation*}
56(8\alpha(2p^2+7v_0^2)-28\beta v_0(f_0+3p)-f_0)x^2y^2+3136v_0xy^3+784f_0y^4
\end{equation*}
which must be searched for a solution.

\section{Numerical Results}
Based on the formulae above, I constructed a very small and straightforward Pari-gp code \cite{PARI2}.

The paper by Komato, Watanabe and Wada consider $p=6079$ with $f_0=1737$ and $v_0=1921$. We have $\alpha=5064$ and $\beta=3761$,
giving  the quartic
\begin{equation*}
d^2=1737x^4+215152x^3y+4376904x^2y^2+6024256xy^3+1361808y^4
\end{equation*}
which, as they state, has a solution when $x=14609, y=-4338$, giving $d=\pm 89697426147$ and $c=313487951$, and hence to the point
on $Y^2=X^3-(7*6079)^2X$ with
\begin{equation*}
X=\frac{-6079 \times 3724108191962548382965^2}{1831552406584343522639^2}
\end{equation*}

From Pari we find that the point has height $57.11$.
The curve has rank $2$, and Magma's {\bf MordellWeilShaInformation} subroutine finds the second generator, with height $56.26$.

We can find much bigger points. For $N=7*971$ we find
\begin{equation*}
X=\frac{-971 \times 22767881370566679951825265^2}{13139349641439531363750001^2}
\end{equation*}
which is a generator of a rank one curve and has height $117.45$. Testing over the first $1000$ primes $p$, the largest point found is at
$N=33509$, with height $122$.

In both cases, the Pari subroutine {\bf hyperellratpoints} proved invaluable, since searching the quartic is easily the most
time-consuming aspect of the computation.

\end{document}